\magnification =1200\hoffset=6mm\voffset=10mm

%\magnification=1000\hoffset=15mm\voffset=15mm

\normalbaselines
\topskip=18pt
\hsize=120mm \vsize=180mm \parindent=5mm \vskip 15mm

\overfullrule=0mm
\def\det{{\rm det}}
\def\d¡{{\rm d¡}}
\def\tr{{\rm tr}}

\def\adots{\mathinner{\mkern2mu\raise1pt\hbox{.}
\mkern3mu\raise4pt\hbox{.}\mkern1mu\raise7pt\hbox{.}}}

\def\GL{{\rm GL }}

\def\L{{\rm L }}

\overfullrule=0mm
\def\det{{\rm det}}
\def\d¡{{\rm d¡}}
\def\tr{{\rm tr}}

\def\adots{\mathinner{\mkern2mu\raise1pt\hbox{.}
\mkern3mu\raise4pt\hbox{.}\mkern1mu\raise7pt\hbox{.}}}

\def\GL{{\rm GL }}

\def\Aut{{\rm Aut}}

\def\L{{\rm L }}

\def\SL{{\rm SL }}

\input AMSsym.def
\input AMSsym.tex
\def\tr{{\rm tr}}

\def\Sym{{\rm Sym}}
\def\n{{\rm n }}
\def\Hom{{\rm Hom}}
\def\dim{{\rm dim}}

\input xy.tex
\xyoption{all}

\input xy.tex
\xyoption{all}

% \the\day /
%\the\month /
%\the\year .%

\bigskip
\bigskip

\centerline{\bf  La th\'eorie des invariants des formes }

\centerline{\bf quadratiques ternaires revisit\'ee.}

\bigskip

\centerline{\bf Bruno Blind*}

\bigskip

\bigskip  

\bigskip

\bigskip

\footnote{}{*Institut Elie Cartan, Nancy-Universit\'e, CNRS, INRIA, Boulevard des Aiguillettes, B.P. 239, F- 54506 Vandoeuvre-l\`es-Nancy, France; \hphantom {}blind@iecn.u-nancy.fr}

\rightline {{\it L'histoire qui va suivre est v\'eridique dans ses moindres d\'etails,}}

\rightline {{\it \`a moins qu'une affreuse erreur n'ait tout fauss\'e depuis le d\'ebut.}}

\rightline {Roland Topor}

\bigskip

\bigskip

\bigskip

\noindent{\bf Abstract:} The simultaneous invariants of $2$, $3$, $4$ and $5$ ternary quadratic forms under the group $\SL(3, {\Bbb C})$ were given by several authors (P. Gordan, C. Ciamberlini, H.W. Turnbull, J.A Todd), utilizing the symbolic method. Using the Jordan algebra structure of the space of ternary quadratic forms, we give these invariants explicitely.

\bigskip

\bigskip

\bigskip

\noindent{\bf R\'esum\'e:} Gr\^ace \`a la structure naturelle d'alg\`ebre de Jordan que porte l'espace $V = \Sym(3,{\Bbb C})$, nous identifions  les invariants, sous $\SL(3, {\Bbb C})$, de plusieurs formes quadratiques ternaires trouv\'es \`a l'aide de la m\'ethode symbolique par divers auteurs (P. Gordan, C. Ciamberlini, H.W. Turnbull, J.A Todd et d'autres).

\bigskip

\bigskip

\bigskip

\noindent{\bf Classification AMS:} 13A50; 15A72; 17C99

\bigskip

\bigskip

\noindent{\bf Mots clefs:} Th\'eorie classique des invariants, Formes quadratiques ternaires; Alg\`ebre de Jordan.

\vfill \eject 

\bigskip

\beginsection {1. Introduction.}

La d\'etermination des invariants (et plus g\'en\'eralement des covariants) de plusieurs formes quadratiques ternaires sous l'action du groupe $G = \SL(3, {\Bbb C})$ est un des probl\`emes 
abondamment  \'etudi\'es par la th\'eorie classique des invariants. Le cas des deux formes est trait\'e par P. Gordan [Gor.1],  il montre en particulier que $A(2)$ (nous notons $A(p)$ l'alg\`ebre des invariants de $p$ formes) est engendr\'ee par $4$ invariants. Le cas $p=3$ est abord\'e en 1886 par C. Ciamberlini (et d'autres), qui trouve que $A(3)$ est engendr\'ee par $11$ invariants [Cia.]. Les r\'esultats de ces deux auteurs ont \'et\'e pr\'ecis\'es en 1923 par Van der Waerden [Wae.]. 
En 1948,  J.A Todd, dans deux longs articles [Tod.1 et Tod.2], examine et compl\`ete le travail effectu\'e en 1910 par H.W. Turnbull [Tur.] sur les cas restants ($p=4$ et $p=5$). Tous ces auteurs  donnent les invariants sous forme symbolique, et, \`a ma connaissance, leur forme explicite, en particulier pour $A(4)$ et $A(5)$, n'a pas \'et\'e donn\'ee.

Un des buts de ce travail est de donner ces invariants sous forme explicite.  A cette fin, nous utilisons la structure d'alg\`ebre de Jordan naturelle que porte l'espace $V$ des formes quadratiques ternaires, cet espace $V$ apparait en effet dans la s\'erie des quatre alg\`ebres de Jordan simples de rang trois, s\'erie qui comporte l'alg\`ebre d'Albert. Pour chacune de ces alg\`ebres, nous pouvons d\'efinir un analogue naturel pour le groupe $G$ et  \'enoncer un probl\`eme  d'invariants, pour lequel les polyn\^omes  trouv\'es dans le cas des formes quadratiques ternaires nous donnent encore des invariants.

Apr\`es les pr\'eliminaires, nous red\'emontrons aux paragraphes 3 et 5  les r\'esultats de P. Gordan et C. Ciamberlini. Le cas $p=2$  est trait\'e dans le paragraphe 3 et est une belle application de la m\'ethode g\'eom\'etrique initi\'ee par T. Vust [Vus.].  Au paragraphe 4, gr\^ace \`a une \'etude du nilc\^one de $pV$, nous obtenons suffisamment d'informations sur un syst\`eme de param\`etres homog\`enes pour en d\'eduire explicitement  
la s\'erie de Poincar\'e de l'alg\`ebre $A(3)$. Nous exploitons au paragraphe 5 cette connaissance de la s\'erie de Poincar\'e pour red\'emontrer les r\'esultats de C. Ciamberlini concernant le cas des trois formes. L'alg\`ebre $A(3)$ est engendr\'ee par $11$ invariants soumis \`a une relation cubique, nous  indiquons en appendice une m\'ethode pour \'ecrire explicitement cette relation.

Nous n'arrivons malheureusement pas \`a  red\'emontrer les r\'esultats de H.W. Turnbull et J.A. Todd concernant $A(4)$ et $A(5)$. Aussi, dans le paragraphe 6, pour le cas des quatre formes, nous utilisons le fait que $A(4)$ est engendr\'ee par ses \'el\'ements de degr\'e inf\'erieur \`a $6$ pour pouvoir les \'ecrire  explicitement, et dans le cas des cinq formes, nous identifions (nous avons mis les d\'etails en appendice) un invariant donn\'e sous forme symbolique par Todd pour conclure. 

Tout au long de ce travail, nous faisons usage du logiciel LIE d\'evelopp\'e par Marc A.A. van Leeuwen, Arjeh M. Cohen et Bert Lisser, nous le d\'esignerons simplement par LIE dans la suite.

Je remercie F. Chargois, J-L. Clerc et P. Y. Gaillard qui m'ont encourag\'e et avec qui j'ai eu des discussions toujours stimulantes. Mes remerciements vont  \'egalement \`a W. Bertram, O. Hijazi et G. Rousseau pour leurs critiques et encouragements.

\bigskip

\bigskip

\bigskip

\bigskip

\beginsection {2. Pr\'eliminaires.}

\bigskip

\noindent {\it 2.1. Les invariants simultan\'es des formes quadratiques ternaires.}

Nous faisons op\'erer le groupe $G = \SL(3, {\Bbb C})$ sur l'espace $V = \Sym(3,{\Bbb C})$ des matrices $(3,3)$ sym\'etriques (identifi\'ees aux formes quadratiques ternaires) par l'action:
$$g.x = gxg^t$$
cette repr\'esentation, not\'ee $\phi_1^2$ dans [Sch.1]), donne lieu \`a une action de $G$ sur l'alg\`ebre ${\Bbb C}[pV]$ des polyn\^omes sur $pV=V\oplus \cdots \oplus V$, et on notera par $A(p)$ la sous-alg\`ebre ${\Bbb C}[pV]^G$ des invariants. C'est une alg\`ebre gradu\'ee par le degr\'e des polyn\^omes, et vu l'action du centre de $G$, il est facile de voir que le degr\'e d'un invariant est  un multiple de $3$. On a donc la d\'ecomposition:
$$A(p) = \bigoplus _{n\geq 0}A_{3n}(p),$$
$A_{3n}(p)$ \'etant le sous espace vectoriel de $A(p)$ constitu\'e par les invariants homog\`enes de degr\'e $3n$.

Il est classique de voir $pV$ comme un 
$\GL(p, {\Bbb C}) \times G$-module, d'o\`u  l'identification de 
${\Bbb C}[pV]$ avec le $\GL(p, {\Bbb C}) \times G$-module $S^{\bullet} (\psi_1(p) \otimes V^{*})$, $\psi_1(p)$ d\'esignant la repr\'esentation standard de $\GL(p, {\Bbb C})$ (on reprend ici les notations de [Sch.2] page 192, disons bri\`evement que  $\psi_i(p) = \Lambda^i(\psi_1(p)$ pour $i\geq 0$, et que $\psi_{(a_1, \cdots , a_m, 0)} $ est la composante de plus haut poids dans le produit tensoriel $S^{a_1}(\psi_1(p)) \otimes \cdots \otimes S^{a_m}(\psi_m(p))$). Le th\'eor\`eme de Cauchy en donne la d\'ecomposition en modules irr\'eductibles, on peut ainsi obtenir la d\'ecomposition du  $\GL(p, {\Bbb C})$-module $A(p)$, on a:
$$A(p)  = \bigoplus_{(a)=(a_1, \cdots ,a_6, 0) } \psi_{(a)}(p) \otimes  \psi_{(a)}(V^{*})^{G}$$
(voir [Sch.2], page 196 ou [Pro.], page 386 pour tout ceci). Une des cons\'equences de cette d\'ecomposition est que pour $p>6$, $A(p)$ est engendr\'ee par les polaris\'es des \'el\'ements de $A(6)$ ([Pro.], page 386 th\'eor\`eme 1); en fait ici il nous suffit de connaitre $A(5)$ ([Pro.], page 386 th\'eor\`eme 2).

  La graduation de $A(p)$ est pr\'eserv\'ee par l'action de $\GL(p, {\Bbb C})$, en particulier on a:
$$A_{3}(p) = \bigoplus_{\scriptstyle (a) = (a_1, a_2, a_3, 0, \cdots )\atop \scriptstyle a_1+2a_2+3a_3 = 3} \psi_{(a)}(p) \otimes  \psi_{(a)}(V^{*})^{G}$$

ici $(a) = (3, \cdots , 0)$, ou  $(a) = (1,1, \cdots, 0)$, ou $(a) = (0,0,1, \cdots, 0)$ les deux derniers cas \'etant exclus par LIE par exemple, et l'on obtient pour le $\GL(p, {\Bbb C})$-module $A_3(p)$:
$$A_3(p) = \psi_{(3, \cdots , 0)}(p) = \Sym^3({\Bbb C}^p).$$

On sait que  la repr\'esentation  $2\phi_1^2$ est cor\'eguli\`ere mais que $3\phi_1^2$ ne l'est pas (cf. [Sch.1]) c'est-\`a-dire que $A(p)$ n'est une alg\`ebre de polyn\^omes que pour $p<3$. Divers auteurs ont donn\'e un syst\`eme g\'en\'erateur de ces alg\`ebres $A(p)$.

Par ailleurs le r\'esultat de P. Gordan dit que l'alg\`ebre $A(2)$ est engendr\'ee par ses \'el\'ements de degr\'e $3$; celui de C. Ciamberlini que l'alg\`ebre $A(3)$ est engendr\'ee par $A_{3}(3) \bigoplus A_{6}(3)$; il en est de m\^eme pour $A(4)$ d'apr\`es H.W.Turnbull et J.A.Todd, tandis que $A(5)$, elle, est engendr\'ee par $A_{3}(5) \bigoplus A_{6}(5) \bigoplus A_{9}(5)$. Tous ces auteurs utilisent la m\'ethode symbolique \'elabor\'ee au XIX si\`ecle.

\bigskip

\noindent {\it 2.2. La m\'ethode symbolique.}

Nous d\'ecrivons ici tr\`es sommairement la notation symbolique employ\'ee par les auteurs cit\'es plus haut pour \'enum\'erer les invariants. On trouvera dans [Gra.-You.] une exposition classique (voir aussi [Gor.2]); pour une pr\'esentation moderne et rigoureuse, nous renvoyons le lecteur \`a [Gro.- Rot.- Ste.] (voir aussi l'article de synth\`ese [Rot.- Stu.]).

Aux \'el\'ements $x, y, z, \cdots $ de $V$ sont associ\'es (en nombre arbitraire) des vecteurs symboliques
$a = (a_1, a_2, a_3)$, $b = (b_1, b_2, b_3)$ etc. qui repr\'esentent $x$;  $a' = (a'_1, a'_2, a'_3)$, $b' = (b'_1, b'_2, b'_3)$ etc. qui repr\'esentent $y$; $a'' = (a''_1, a''_2, a''_3)$, $b'' = (b''_1, b''_2, b''_3)$ etc. qui repr\'esentent $z$... On forme l'alg\`ebre de polyn\^omes 

\noindent ${\Bbb C}[\{a_i\}, \{a'_i\}, \{a''_i\}, \cdots , \{b_i\}, \{b'_i\}, \{b''_i\}, \cdots]$ appel\'e espace ombral et on d\'efinit l'application ombrale ${\cal U}$:
$${\cal U}: {\Bbb C}[\{a_i\}, \{a'_i\}, \{a''_i\}, \cdots , \{b_i\}, \{b'_i\}, \{b''_i\}, \cdots] \mapsto {\Bbb C}[pV]$$
c'est une application lin\'eaire dont nous rappelons seulement  les propri\'et\'es suivantes: si $x$ (resp $y$) est la matrice $(x_{ij})$ (resp $(y_{ij})$), on a:
$${\cal U}(a_i^2) = {\cal U}(b_i^2) = \cdots = x_{ii}; \   \  {\cal U}(a_ia_j) = {\cal U}(b_ib_j) = \cdots = x_{ij} (i\not = j)$$
$${\cal U}(a^r_ia^s_j) = {\cal U}(b^r_ib^s_j) = 0 \ \ \hbox {si}  \ \  r+s \not = 2$$
$${\cal U}((a'_i)^2) = {\cal U}((b'_i)^2) = \cdots = y_{ii}; \   \  {\cal U}(a'_ia'_j) = {\cal U}(b'_ib'_j) = \cdots = y_{ij} (i\not = j)$$
(voir [Rot.- Stu.] page 8 pour une d\'efinition compl\`ete de ${\cal U}$). D'autre part, si $u$, $v$, $w$ sont trois vecteurs symboliques, leur d\'eterminant, appel\'e crochet, sera not\'e $[u,v,w]$. Le point clef de la m\'ethode symbolique est alors le th\'eor\`eme suivant (voir les r\'ef\'erences cit\'ees plus haut):

\proclaim Th\'eor\`eme 2.1. Tout \'el\'ement de $A(p)$ est de la forme ${\cal U}(Q)$ o\`u $Q$ est un polyn\^ome en les crochets.

\noindent On v\'erifie par exemple que si les trois vecteurs symboliques repr\'esentent l'\'el\'ement $x$, leur crochet redonne essentiellement l'invariant $\det x$, d'une mani\`ere plus pr\'ecise:

$${\cal U}([a,b,c]^2) = {1\over 6}\det x .$$

\noindent Introduisons, pour finir, une derni\`ere notation symbolique: \`a partir de deux vecteurs symboliques 
$u = (u_1, u_2, u_3)$ et $v = (v_1, v_2, v_3)$, on construit un
troisi\`eme triplet symbolique (produit vectoriel de $u$ avec $v$) de composantes:
$$ {u_2}{v_3} - {u_3}{v_2}, \   \   \  
{u_3}{v_1} - {u_1}{v_3}, \  \  \
 {u_1}{v_2} -
{u_2}{v_1}$$
quand $u$ et $v$ repr\'esentent  $x$, on notera cet \'el\'ement par $\alpha$, et on \'ecrira le crochet $[a,b,c]$ comme produit scalaire de $a$ avec ${\alpha}$: $[a,b,c] = a_{\alpha}$. P. Gordan montre qu'un syst\`eme g\'en\'erateur minimal de l'alg\`ebre $A(2)$ s'\'ecrit dans ces notations:
$$a_{\alpha}^2,\  \  \  (a'_{\alpha'})^2,\  \  \  (a'_{\alpha})^2,\  \  \  a_{\alpha'}^2 .$$
Les r\'esultats de C. Ciamberlini nous donnent un syst\`eme g\'en\'erateur minimal de l'alg\`ebre $A(3)$  form\'e de onze invariants:
$$a_{\alpha}^2, a_{\alpha'}^2 , a_{\alpha''}^2, (a'_{\alpha'})^2, (a'_{\alpha})^2, (a'_{\alpha''})^2, (a''_{\alpha''})^2, (a''_{\alpha})^2, (a''_{\alpha'})^2, [a,a',a'']^2, [\alpha, \alpha',\alpha'']^2$$

\bigskip

\noindent {\it 2.3. $V$ vue comme alg\`ebre de Jordan.}

On s'aper\c coit que ces invariants donn\'es sous forme symbolique peuvent s'exprimer sous forme concr\`ete en utilisant la structure d'alg\`ebre de Jordan de  l'espace  $V = \Sym(3,{\Bbb C})$ muni du produit:
$$x \bullet y = {1\over 2}(xy+yx).$$
De plus, les formules obtenues ont un sens dans la situation plus g\'en\'erale des alg\`ebres de Jordan sur ${\Bbb C}$ simples et de rang $3$. Nous d\'etaillons bri\`evement cette situation, en renvoyant le lecteur \`a [Far.-Kor.] et \`a [Spr.] pour les d\'efinitions et les principaux faits  concernant ces alg\`ebres.

Soit donc $V$ une alg\`ebre de Jordan sur ${\Bbb C}$ simple, de rang $3$ et de dimension $n$; $V$ est munie de la forme bilin\'eaire d\'efinie positive $<x,y> = {3\over n}\tr L(x\bullet y) $ o\`u $\L(x)$ d\'esigne l'endomorphisme multiplication par $x$: $\L(x): y\mapsto x\bullet y$, cette forme v\'erifie 
$$<x\bullet z, y> = <x, y\bullet z>.$$
On note par $\det$ le d\'eterminant de $V$, c'est un polyn\^ome irr\'eductible et homog\`ene de degr\'e le rang de $V$, c'est \`a dire ici $3$; un \'el\'ement $x$ de $V$ est inversible si $\det x \not = 0$, et dans ce cas l'inverse $x^{-1}$ s'\'ecrit sous la forme
$$x^{-1} = {\n(x)\over \det(x)}$$
o\`u $\n$ est une application quadratique de $V$ dans $V$, on d\'esigne par
$\times$ l'application bilin\'eaire sym\'etrique associ\'ee \`a $\n$: $x\times y = \n(x+y) -\n(x)-\n(y)$. Par polarisation du polyn\^ome $\det$, nous obtenons une forme trilin\'eaire $f$ sur $V$ telle que $f(x,x,x)$ = $6 \ \det x$, et nous avons l'importante relation (voir [Spr.], chapitre 4, formule (4), page 55):

$$f(x, y, z) =<x \times y, z>.$$
Soit  $G$ le sous-groupe  de $\GL(V)$ laissant invariant la fonction $\det$; c'est un sous-groupe du groupe de structure de $V$ (voir [Spr.], Proposition 12.3, page 123) qui contient le groupe $\Aut(V)$ des automorphismes de l'alg\`ebre $V$. La forme trilin\'eaire $f$ est invariante par $G$; d'autre part, il n'est pas difficile de constater que l'application $\n$ v\'erifie la propri\'et\'e de $G$-\'equivariance suivante
$$\n(g.x) = (g^{-1})' \n(x)$$
o\`u $g'$ d\'esigne l'adjoint de $g$ par rapport \`a la forme bilin\'eaire $<, >$ (voir [Far.-Kor.] page 148). A partir de la forme trilin\'eaire $f$ et de l'application bilin\'eaire $\times $, nous pouvons  donc facilement construire des invariants du groupe $G$, par exemple le polyn\^ome $f(\n(x), \n(y),\n(z))$ ainsi que tous ses polaris\'es sont des invariants. Ces invariants ont \'egalement \'et\'e consid\'er\'es par A.V. Iltyakov dans [Ilt.].

Dans le cas de l'espace  $V = \Sym(3,{\Bbb C})$ muni du produit de Jordan
$x \bullet y = {1\over 2}(xy+yx)$, la fonction $\det$ est le  d\'eterminant usuel des matrices et $n(x)$ n'est autre que la transpos\'ee de la matrice des cofacteurs de la matrice $x$ (l'inverse de Jordan et l'inverse usuel coincident); la forme bilin\'eaire $<x, y> $ est la trace habituelle de $xy$ et $G$ est le  groupe $\SL(3, {\Bbb C})$ agissant sur $V$ par $g.x = gxg^t$. 

Le lecteur trouvera dans [Far.-Kor.] et dans [Spr.] la classification (il y en a $4$) des alg\`ebres de Jordan sur ${\Bbb C}$ simples et de rang $3$, rappelons seulement qu'outre le cas des formes quadratiques ternaires que nous envisageons dans ce travail, il y a l'alg\`ebre d'Albert et que dans ce cas le groupe $G$ est le groupe complexe simplement connexe $E_6$.

\bigskip

\bigskip

\bigskip

\bigskip

\beginsection {3. L'alg\`ebre $A(2)$.}

\bigskip

On se propose dans ce paragraphe de red\'emontrer le th\'eor\`eme de P. Gordan concernant $A(2)$ (cf [Gor.1]); nous en donnons une preuve qui se g\'en\'eralise imm\'ediatement au cas des alg\`ebres de Jordan sur ${\Bbb C}$ simples et de rang $3$ \'evoqu\'ees plus haut. Dans toute la suite, on fixe un rep\`ere de Jordan $\{e_1, e_2, e_3\}$ de $V$.
\proclaim Th\'eor\`eme 3.1. L'alg\`ebre $A(2)$ est une alg\`ebre de polyn\^omes engendr\'ee par les quatres polyn\^omes:
$$\det \ x, \   \   \det \ y, \  \  f(x,x,y), \  \  f(x,y,y).$$

La d\'emonstration consiste \`a utiliser la m\'ethode g\'eom\'etrique de T. Vust (cf [Vus.]). On identifie 
$2V$ \`a $\Hom ({\Bbb C}^2, V)$: l'\'el\'ement $(x,y)$ de $2V$ est vu comme l'application lin\'eaire $\alpha_{(x,y)}$ d\'efinie par:
$\alpha_{(x,y)}(a,b) = ax+by$,
 et on consid\`ere l'application $\chi $ de $2V$ dans l'espace $\Sym^3({\Bbb C}^2)$ des polyn\^omes homog\`enes de degr\'e $3$ sur ${\Bbb C}^2$ donn\'ee par la formule:
$\chi (x,y) = \det \circ \alpha_{(x,y)} $,
c'est-\`a-dire que $\chi (x,y)$ est le polyn\^ome \`a deux variables $a$ et $b$ suivant:
$$(a,b) \mapsto \det (ax+by) = a^3 \det x + b^3 \det y +{a^2b\over 2} f(x,x,y) + {ab^2\over 2} f(x,y,y).$$
Il est alors clair que le morphisme  alg\'ebrique $\chi $ est $G$-invariant. Il se factorise donc suivant le diagramme:

$$\xymatrix{
&2V\ar[ld]_{\pi}\ar[rd]^{\chi}&\\
2V //G\ar@{-->}[rr]_{\Psi}&&{\hbox{Sym}}^3({\Bbb C}^2)}$$
o\`u l'espace $2V //G$ est la vari\'et\'e alg\'ebrique affine irr\'eductible associ\'ee \`a l'alg\`ebre (int\`egre, de type fini) $A(2)$ et l'application $ \pi: 2V \rightarrow 2V //G$ est le morphisme canonique. 

La strat\'egie est alors la suivante: on montre que $\Psi$ est surjective et birationnelle; il en r\'esultera que c'est un isomorphisme alg\'ebrique car $\Sym^3({\Bbb C}^2)$ est une vari\'et\'e normale (voir par exemple [Bri.], lemme 1 page 132), et, vu la forme de l'application $\chi $, le th\'eor\`eme 3.1 s'en suivra.

Il n'est pas difficile de montrer que $\chi$, et par cons\'equent $\Psi$, est surjective: en effet, consid\'erons  un \'el\'ement $g$ de $\Sym^3({\Bbb C}^2)$, comme on a $g(a,b) = b^3g({a\over b},1)$, le polyn\^ome $g$ peut s'\'ecrire sous la forme:
$$g(a,b) = \lambda b^{3-N}\prod_1^N (a-\alpha_ib)$$
o\`u $\lambda$ est un scalaire non nul, et avec $0\leq N \leq 3$ (un produit vide \'etant \'egal \`a $1$). En posant alors $x = \sum_1^N e_i$ (on prend $x=0$ si $N=0$) et $y = -\sum_1^N\alpha_ie_i + e_{N+1} \cdots +e_3$, il vient 
$$\det \circ \alpha_{(x,y)}  = {1\over \lambda}f$$
et donc $gÊ\in Im \chi$, puisque cette image  est ${\Bbb C}^*$ invariant.

Le morphisme alg\'ebrique $\Psi$ \'etant surjectif, il est dominant; pour montrer qu'il est birationnel il nous suffit (voir par exemple [Vin.], lemme 1 page 252) d'exhiber un sous-ensemble $\cal V$ dense de $\Sym^3({\Bbb C}^2)$ tel que pour tout $v$ de $\cal V$ la fibre $\Psi^{-1}(v)$ soit r\'eduite \`a un point. On prend alors
$$\cal V = \chi (\cal U)$$
$\cal U$ \'etant le sous-ensemble dense de $2V$ constitu\'e des couples $(x,y)$ avec $x$ inversible et $y$ poss\'edant trois valeurs propres distinctes. Soit $(x_0, y_0)$ un \'el\'ement de $\cal U$, on va montrer que si $(x, y) \in 2V$ est tel que 
$$\chi (x_0,y_0) = \chi (x,y) \leqno (*) $$
alors $(x, y)$ et $(x_0, y_0)$ sont dans la m\^eme $G$-orbite, ce qui entrainera bien que la fibre $\Psi^{-1}(\chi (x_0,y_0))$ est r\'eduite \`a $\pi (x_0,y_0)$. Tout d'abord, on peut supposer que $\det x_0 = 1$ sans  restreindre  la g\'en\'eralit\'e du probl\`eme; on peut alors (voir [Far.-Kor.], proposition VIII.3.5, page 153), en se d\'epla\c cant dans la $G$-orbite de $(x_0, y_0)$, supposer que $x_0 = e$. La relation $  (*) $ entraine alors que $\det x = 1$, on peut donc encore, en se d\'epla\c cant cette fois ci dans la $G$-orbite de $(x, y)$ (et quitte \`a changer de notation) remplacer  $(x, y)$ par $(e, y)$; la relation $  (*) $ nous donne en particulier pour tout complexe $a$:
$$\det \ (ae+y_0) =  \det \ (ae+y)$$
et ceci implique que $y$ poss\`ede les m\^emes  valeurs propres distinctes $\lambda_i, i=1, 2, 3$ que $y_0$. On sait alors (voir [Far.-Kor.], proposition VIII.3.2, page 151) que les \'el\'ements $y_0$ et $y$ peuvent s'\'ecrire sous la forme:
$$y_0 =  \lambda_1c_1 +\lambda_2c_2 + \lambda_3c_3, \  \  y =  \lambda_1c'_1 +\lambda_2c'_2 + \lambda_3c'_3$$
o\`u $(c_1, c_2, c_3)$ et  $(c'_1, c'_2, c'_3)$ sont deux rep\`eres de Jordan de $V$, mais l'on sait (voir 
[Far.-Kor.], th\'eor\`eme IV.2.5, page 71) que le groupe $\Aut(V)$ des automorphismes de $V$ op\`ere transitivement sur les rep\`eres de Jordan. Il r\'esulte de tout ceci que $(x, y)$ et $(x_0, y_0)$ sont bien dans la m\^eme $G$-orbite.

\bigskip

\bigskip

Le r\'esultat de P. Gordan nous permet de trouver  la dimension de Krull $d(p)$ de $A(p)$ pour $p\geq 2$; en effet comme $d(2) = 4$, il r\'esulte de la formule:
$$d(2) = \dim 2V - \max \dim (G.v)$$
que la dimension maximale d'une $G$-orbite $G.v$ dans $2V$ est  $8$, ce qui est la dimension du groupe $G$; par cons\'equent la dimension maximale d'une $G$-orbite dans $pV$ est \'egalement de $8$ et donc on a la:

\proclaim Proposition 3.3. La dimension de Krull $d(p)$ de $A(p)$ ($p\geq 2$) est:
$$d(p) = 6p-8.$$

\bigskip

\bigskip

Pour terminer ce paragraphe, donnons la structure des $\GL(2,{\Bbb C})$ modules $A_3(2)$, $A_6(2)$ et $A_9(2)$. Nous avons (voir le paragraphe 2):

$$A_3(2)  = \Sym^3({\Bbb C}^2) = \psi_{(3, \cdots , 0)}(2).$$
Par LIE nous obtenons:
$$A_6(2) = {\hbox{Sym}}^2({\hbox{Sym}}^3({\Bbb C}^2)) = \psi_{(6,0,  \cdots , 0)}(2) \oplus \psi_{(2,2,0, \cdots , 0)}(2)$$
$$A_9(2) = {\hbox{Sym}}^3({\hbox{Sym}}^3({\Bbb C}^2)) = \psi_{(9,0,  \cdots , 0)}(2) \oplus \psi_{(5,2,0, \cdots , 0)}(2)\oplus \psi_{(3,3,0, \cdots , 0)}(2)\oplus \psi_{(0,3,1, \cdots , 0)}(2)$$

\bigskip

\bigskip

\bigskip

\bigskip

\beginsection {4. Syst\`eme de param\`etres homog\`enes.}

\bigskip

 Notons  $B(p)$ la sous-alg\`ebre de $A(p)$  engendr\'ee par les polynomes $f(x_i, x_j, x_k), $
 
\noindent $ 1\leq i,j,k\leq p $. Le but de ce paragraphe est de montrer le th\'eor\`eme:

\proclaim Th\'eor\`eme 4.1. L'alg\`ebre $A(p)$ est enti\`ere sur $B(p)$.

Il est classique (voir par exemple [Der.-Kem.], lemme 2.4.5, page 60) que ce th\'eor\`eme r\'esulte du lemme suivant (o\`u ${\cal N}_{pV}$ ($p\geq 1$) d\'esigne le nilc\^one de $pV$):

\proclaim Lemme 4.2. On a:
$$ {\cal N}_{pV} = \{ (x_1, \ldots, x_p) \in pV \ \ / f(x_i, x_j, x_k) = 0 \ \ 1\leq
i,j,k\leq p \}.$$

On peut bien entendu supposer, sans restreindre la g\'en\'eralit\'e, que $p\geq 3$. Soit $\{x_1,x_2, \cdots ,x_p\}$ une famille satisfaisant aux relations $f(x_i,x_j,x_k) = 0$, par application du  crit\`ere de Hilbert-Mumford (voir [Bri.] page 139, ou bien [Der.-Kem.] page 60) il suffit de montrer qu'il existe un sous-groupe \`a un param\`etre $\lambda$ de $G$ tel que 
$$\lim_{t\rightarrow 0} \lambda (t) x_i = 0, \   \   i = 1, 2, \cdots ,p.$$
A conjugaison pr\`es, les sous-groupes \`a un param\`etre de $G$ sont les $\lambda_{(n_1, n_2, n_3)}$:

$$\lambda_{(n_1, n_2, n_3)}: t\in {\Bbb C}^* \mapsto (\pmatrix{t^{n_1}&0&0\cr
0&t^{n_2}&0\cr
0&0&t^{n_3} \cr}, \   \   n_3 = -n_1-n_2$$
op\'erant sur $V$ par:
$$\lambda_{(n_1, n_2, n_3)}(t). (x_{ij}) = (t^{n_i+n_j} x_{ij}).$$
La d\'emonstration consistera alors \`a montrer que dans la $G$-orbite de $\{x_1,x_2, \cdots ,x_p\}$ il existe soit une famille form\'ee de matrices toutes de la forme:
$$\pmatrix{\star &\star & 0 \cr
 \star &\star & 0 \cr
 0 & 0 & 0 \cr},$$
soit une famille form\'ee de matrices toutes de la forme:

$$\pmatrix{\star &\star & \star \cr
 \star &0 & 0 \cr
 \star & 0 & 0 \cr}.$$
Pour ces deux formes, il n'est pas difficile d'exhiber explicitement un  sous-groupe \`a un param\`etre $\lambda_{(n_1, n_2, n_3)}$ qui convient.

\noindent Nous avons choisi de donner ici  une preuve qui se g\'en\'eralise \`a la situation des alg\'ebres de Jordan simples et de rang trois, et qui repose sur l'emploi des transformations de Frobenius (voir [Far.-Kor.] page 106). On \'ecrira indiff\'eremment un \'el\'ement $y$ de $V$ sous la forme matricielle $y = (y_{ij})$ ou dans la d\'ecomposition de Peirce  $y = y_{11}e_1 + y_{22}e_2 +y_{33}e_3+{\bf y}_{12} + {\bf y}_{13} + {\bf y}_{23}$ associ\'ee au rep\`ere de Jordan $(e_1, e_2, e_3)$ (avec $e_1 = \pmatrix{1&0 & 0 \cr
0 &0 & 0 \cr
0 & 0 & 0 \cr}$, $ {\bf y}_{12} = \pmatrix{0 &y_{12} & 0 \cr
 y_{12} &0 & 0 \cr
0 & 0 & 0 \cr}$, etc). On peut, bien entendu, donner  une preuve plus \'el\'ementaire de ce qui va suivre. 

\bigskip

$\bullet$  On suppose dans un premier temps que l'un des $x_i$, par exemple $x_1$, soit de rang $2$; par action du groupe $G$, et par multiplication d'un scalaire, on peut se ramener au cas o\`u $x_1 = e_1+ e_2$. Soit alors $y$ un autre \'el\'ement $x_i$ de la famille; comme $\n(e_1+e_2) = e_3$,
la relation $f(x_1, x_1, y) = 0$ donne
$$y_{33} = 0$$
 et en calculant alors les composantes de $\n(y)$, la relation $f(x_1, y, y) = 0$ entraine que:
$$y_{23}^2 = - y_{13}^2.$$
Si l'on a $ y_{23} = 0$ pour chacun des  $y$, tous les  \'el\'ements de la famille sont de la forme:
$$\pmatrix{\star &\star & 0 \cr
 \star &\star & 0 \cr
 0 & 0 & 0 \cr}$$
et nous avons fini; on peut donc supposer qu'il existe un \'el\'ement $y = (y_{ij})$ de la famille avec $y_{23} \not = 0$; nous allons montrer que dans ce cas nous pouvons supposer que la famille soit du type $(x, y', \cdots )$ avec:
 
 $$x = \pmatrix{1&a&0 \cr
a&0& 0 \cr
0 & 0 &0\cr}, \ \     y' = \pmatrix{*&0&a'\cr
0 & 0 & 0 \cr
a'&0&0 \cr}$$
$a$ et $a'$ \'etant deux nombres complexes non nuls. Pour ce faire, appliquons  \`a notre famille la transformation de Frobenius $\tau(u)$ avec 
$$u = \pmatrix{0&-{y_{23} \over  y_{13}}&0 \cr
-{ y_{23} \over  y_{13}}&0& 0 \cr
0 & 0 &0\cr}$$
il vient (voir [Far.-Kor.], lemme VI.3.1, page 106):
$$\tau(u)(e_1+e_2) =     e_2 + e_1 + 2 ue_1 + 2 (e_2+e_3)(u(ue_1)) = $$
$$e_2 + e_1 + u + (e_2+e_3)(u^2) = e_1 + u$$
de m\^eme, pour  $y = y_{11}e_1 + y_{22}e_2 +{\bf y}_{12} + {\bf y}_{13} + {\bf y}_{23}$, il vient:
$$\tau(u)(y) = y_{11} \tau (u)(e_1) + \tau(u)({\bf y}_{12} + {\bf y}_{13})+ y_{22}e_2 + {\bf y}_{23} = $$
$$ y_{11}(e_1 + u - e_2) + {\bf y}_{12} + {\bf y}_{13} + 2(e_2+e_3)(u({\bf y}_{12} + {\bf y}_{13})) + y_{22}e_2 +{\bf y}_{23} = $$
$$y_{11}e_1+(y_{22}- y_{11}-2{y_{12} y_{23} \over  y_{13}})e_2   + y_{11}u + {\bf y}_{12} + {\bf y}_{13} + {\bf y}_{23}  + 2u{\bf y}_{13}$$
mais $y_{22}- y_{11}-2{y_{12} y_{23} \over  y_{13}} = 0$ puisque $\det \ y =0$ et $ {\bf y}_{23} + 2u{\bf y}_{13} = 0$, et par cons\'equent:
$$\tau(u)(y) =  y_{11}e_1 + y_{11}u + {\bf y}_{12} + {\bf y}_{13}$$
maintenant, quitte \`a rajouter \`a notre famille une combinaison lin\'eaire de $\tau(u)(y)$ et de $\tau(u)(e_1+e_2)$, nous pouvons supposer que la famille poss\`ede un \'el\'ement tel que $y'$. 

Soit maintenant $z$ un troisi\`eme \'el\'ement de la famille et regardons les  \'equations v\'erifi\'ees par le triplet $(x,y',z)$. De l'\'equation $f(x, x, z) = 0$ nous tirons \`a nouveau $z_{33} = 0$, de $f(y', y', z) =0$ nous pouvons conclure que $z_{22} = 0$, enfin $f(x, y', z) = 0$ nous donne:
$${\bf z}_{23} = 0$$
(en effet, un  calcul simple nous montre que $x\times y' = 2 {\bf x}_{12} {\bf y'}_{13} = \pmatrix{0&0&0 \cr
0 & 0 & aa' \cr
0 &aa'&0 \cr}$).

\bigskip

$\bullet$  Il reste le cas o\`u tous les $x_i$ sont de rang $1$. Par action du groupe $G$, et par multiplication d'un scalaire, on prend $x_1 = e_1$. On regarde la composante $(\ )_{33}$ des \'el\'ements de la famille.

$\mapsto$ Supposons dans un premier temps qu'il existe un $x_i$, avec cette composante non nulle; notons cet \'el\'ement par $y$:

$$y = y_{11}e_1 + y_{22}e_2 + y_{33}e_3 + {\bf y}_{12} + {\bf y}_{13} + {\bf y}_{23}, \  \  y_{33} \not = 0.$$ 
Nous allons, par une transformation de Frobenius $\tau(u)$ laissant fixe $e_1$, changer $y$ en $y' = y_{33}e_3$; en rajoutant ensuite \`a notre famille de d\'epart l'\'el\'ement $e_1 + y'$, nous serons ramen\'e au cas pr\'ec\'edent. Prenons donc $\tau(u)$ avec $u = -{1\over y_{33}}({\bf y}_{13} + {\bf y}_{23})$; cette transformation laisse fixe $e_1$, $e_2$ et les composantes $(\ )_{12}$, par cons\'equent:
$$\tau(u)(y) = y_{11}e_1 + y_{22}e_2 + {\bf y}_{12} + {\bf y}_{13} + {\bf y}_{23}  + 2 \L(e_1+e_2)[u({\bf y}_{13} + {\bf y}_{23})] +$$
$$y_{33}(e_3 +2ue_3) + 2y_{33}\L(e_1+e_2)[u(ue_3)]$$
cela donne apr\`es un petit calcul et compte tenu du choix de $u$:
$$\tau (u)(y) = y' = y'_{11}e_1 + y'_{22}e_2 + y_{33}e_3 + {\bf y}'_{12}.$$
Cet \'el\'ement $y'$ \'etant de rang $1$ puisque les transformations de Frobenius conservent le rang, il vient:
$$y' = y_{33}e_3$$
comme souhait\'e.

\bigskip

$\mapsto$ Supposons enfin que la composante $(\ )_{33}$ de tous les $x_i$ soit nulle. Si pour tous ces 
$x_i$ les composantes $(\ )_{23}$ et $(\ )_{22}$ sont nulles, tous les  \'el\'ements de la famille seront de la forme:
$$\pmatrix{\star &\star & \star \cr
 \star &0 & 0 \cr
 \star & 0 & 0 \cr}$$
et nous avons fini. Soit donc 
$$y = y_{11}e_1 + y_{22}e_2  + {\bf y}_{12} + {\bf y}_{13} + {\bf y}_{23}$$
un \'el\'ement de la famille avec $y_{22}$ et $ y_{23}$ non nuls simultan\'ement. Nous allons, toujours par une transformation de Frobenius $\tau(u)$, ramener $y$ \`a $y'$ avec $y'_{33} \not = 0$. On prend cette fois fois ci $\tau(u)$, avec $u = {\bf u}_{23}$, laissant fixe $ y_{11}e_1 + {\bf y}_{13}$, et qui fait apparaitre dans $\tau(u)(y) = y'$ le terme
$$y'_{33} = 2u_{23} y_{23}  + y_{22}{u_{23}}^2.$$
Vu les hypoth\`eses faites, nous pouvons rendre cette expression non nulle, et ainsi nous ramener au cas pr\'ec\'edent, ce qui ach\`eve la d\'emonstration.

\bigskip

\noindent {\bf Remarque 4.3.} Dans un travail en cours de r\'edaction, nous montrons un analogue du th\'eor\`eme 4.1  pour les autres alg\`ebres de Jordan simples  de rang trois. Dans ces cas,  il faut adjoindre \`a  la sous-alg\`ebre $B(p)$ d'autres invariants (par exemple dans le cas de l'alg\`ebre d'Albert o\`u $G$ est le groupe $E_6$, il faut ajouter un invariant de degr\'e 9),  la preuve du lemme 4.2 devient alors nettement plus p\'enible.

Remarquons enfin pour terminer ce paragraphe, qu'il r\'esulte du th\'eor\`eme 4.1 que l'alg\`ebre $A(p)$ poss\`ede un syst\`eme de param\`etres homog\`enes constitu\'e de $6p-8$ \'el\'ements de degr\'e $3$ qui sont des combinaisons lin\'eaires des polyn\^omes $f(x_i, x_j, x_k), \  1\leq i,j,k\leq p$ (voir par exemple [Der.-Kem.], lemme 2.4.7, page 61) .

\bigskip

\bigskip

\bigskip

\bigskip

\beginsection {5. L'alg\`ebre $A(3)$.}

\bigskip

Dans ce paragraphe, nous  red\'emontrons  le r\'esultat de C. Ciamberlini  concernant l'alg\`ebre $A(3)$. Ici la connaissance de la s\'erie de Poincar\'e  nous suffira pour d\'eterminer l'alg\`ebre. 

Notons par $f_1, \cdots f_{10}$ les dix polyn\^omes:
$$\det x, \   \   \det y, \  \ \det z, \  \  f(x,x,y), \  \ f(x,x,z)$$
$$  f(y,y,x), \  \  f(y,y,z), \  \  f(z,z,x), \  \  f(z,z,y), \  \  f(x,y,z).$$
Par le paragraphe pr\'ec\'edent, nous savons que ces dix polyn\^omes forment un syst\`eme de param\`etres homog\`enes de $A(3)$; cette alg\`ebre  est donc un ${\Bbb C}[f_1, \cdots , f_{10}] $-module libre de rang $r$: il existe  des invariants $g_2, \cdots  g_r$ tels que:
$$A(3) = {\Bbb C}[f_1, \cdots , f_{10}] \bigoplus {\Bbb C}[f_1, \cdots , f_{10}]g_2\bigoplus \cdots {\Bbb C}[f_1, \cdots , f_{10}]g_r.$$
Nous allons d\'eterminer ce rang en calculant  la  s\'erie de Poincar\'e $P(t)$ de $A(3)$. Cette s\'erie  s'\'ecrit:
$$P(t) = {\sum_{i=0}^l a_{3i}t^{3i} \over (1-t^3)^{10}}$$
puisque les 10 polyn\^omes $f_1, \cdots f_{10}$ sont de degr\'e 3. D'autre part, il r\'esulte des tables de F. Knop et P. Littelmann (voir [Kno.- Lit.]) que:

$$3l = 30 - 18 = 12$$
enfin, les anneaux $A(p)$ \'etant de Gorenstein (cf [Pro.], page 562) nous avons $a_0 = a_{12} = 1$ et $a_{12-3i} = a_{3i}$;  finalement:
$$P(t) = {1+a_3t^3 + a_6t^6 + a_3t^9 + t^{12}\over (1-t^3)^{10}}.$$
On a $\dim A_3(3) = \dim  \Sym^3({\Bbb C}^3) = 10$ et, par utilisation du logiciel LIE, on obtient: 
$$\dim A_6(3) =  56.$$
On en tire facilement:
$$P(t) = {1+t^6 + t^{12}\over (1-t^3)^{10}}.$$
En particulier, le ${\Bbb C}[f_1, \cdots , f_{10}] $-module libre $A(3)$ est de rang  $r=3$, et il  existe une base  de la forme $(1, g, g^2)$, o\`u $g$ est invariant de degr\'e $6$. L'alg\`ebre $A(3)$ est donc engendr\'ee par ses \'el\'ements de degr\'e $3$ et $6$.  L'espace vectoriel $ \psi_{(6,0,  \cdots , 0)}(3) \oplus \psi_{(2,2,0, \cdots , 0)}(3)$ est de dimension $55$, on cherche donc par LIE un candidat $\psi_{(a_1, a_2, a_3, \cdots , 0)}(3)$ (de dimension $1$) pouvant intervenir dans la d\'ecomposition du 
$\GL (3, {\Bbb C})$ module $A_6(3)$, c'est \`a dire tel que $\psi_{(a_1, a_2, a_3, \cdots , 0)}(V^*) \not = 0$; on trouve:

$$A_6(3) = \psi_{(6,0,  \cdots , 0)}(3) \oplus \psi_{(2,2,0, \cdots , 0)}(3)\oplus \psi_{(0,0,2, \cdots , 0)}(3)$$
$$= {\hbox{Sym}}^2({\hbox{Sym}}^3({\Bbb C}^3)) \oplus \psi_{(0,0,2, \cdots , 0)}(3) = A_3(3)A_3(3)\oplus \psi_{(0,0,2, \cdots , 0)}(3)$$
Il n'est pas difficile de v\'erifier que l'invariant $f_{11} (x,y,z) = f(\n (x), \n (y),\n (z))$ n'appartient pas \`a $A_3(3)A_3(3)$ et que donc il a une projection non nulle sur l'espace $\psi_{(0,0,2, \cdots , 0)}(3)$ de $A_6(3)$. En effet, supposons que $f_{11}$ soit dans $A_3(3)A_3(3)$, nous aurions une relation du type:
$$f_{11} (x,y,z) = A[f(x,x,y)f(y,z,z) + f(x,x,z)f(z,y,y) + f(y,y,x)f(x,z,z)] + $$
$$B f(x,y,z)^2$$
prenons cette relation pour $x = e_1$ et $y = e_2$, il vient $\n(x) = \n(y) = 0$ et donc $Bf(e_1, e_2, z) = 0$ pour tous les $z$, d'o\`u $B=0$, et par cons\'equent:
$$f_{11} (x,y,z) = A[f(x,x,y)f(y,z,z) + f(x,x,z)f(z,y,y) + f(y,y,x)f(x,z,z)]$$
on fait cette fois ci $x = e_1$, on aurait pour tous les $y$ et $z$: $Af(y,y,e_1)f(e_1,z,z) = 0$, ce qui conduit \`a $A=0$, absurde, puisque l'invariant $f_{11}$ n'est pas nul identiquement. Finalement nous avons:

\proclaim Proposition 5.1. $(1, f_{11}, f_{11}^2)$ est une base du ${\Bbb C}[f_1, \cdots , f_{10}] $-module libre $A(3)$.

\noindent {\bf Remarque 5.2.} Il r\'esulte de la forme de la s\'erie de Poincar\'e que l'invariant $f_{11}$ v\'erifie une relation du type:
$$f_{11}^3 + Q_1(f_1, \cdots , f_{10})f_{11}^2 + Q_2(f_1, \cdots , f_{10})f_{11} + Q_3(f_1, \cdots , f_{10}) = 0$$
nous donnerons en appendice une m\'ethode pour \'ecrire  explicitement cette relation. Disons \'egalement que dans le cas des autres alg\`ebres de Jordan simples  de rang trois, l'alg\`ebre $A(3)$ est une alg\`ebre de polyn\^omes engendr\'ee par les onze \'el\'ements $f_1, \cdots , f_{11}$.

Enfin, il est utile de noter pour le prochain paragraphe la structure  du $\GL (3, {\Bbb C})$ module $A_9(3)$:
$$A_9(3) = \psi_{(9,0,  \cdots , 0)}(3) \oplus \psi_{(5,2,0, \cdots , 0)}(3)\oplus \psi_{(3,3,0, \cdots , 0)}(3)\oplus \psi_{(0,3,1, \cdots , 0)}(3)\oplus 2\psi_{(3,0,2, \cdots , 0)}(3)$$
$$= {\hbox{Sym}}^3({\hbox{Sym}}^3({\Bbb C}^3)) \oplus \psi_{(3,0,2, \cdots , 0)}(3)$$
$A_9(3)$ contient $\psi_{(9,0,  \cdots , 0)}(3) \oplus \psi_{(5,2,0, \cdots , 0)}(3)\oplus \psi_{(3,3,0, \cdots , 0)}(3)\oplus \psi_{(0,3,1, \cdots , 0)}(3)\oplus \psi_{(3,0,2, \cdots , 0)}(3)$;  LIE nous donne leur dimension: $230 $ et $220 $ et permet de v\'erifier que $\psi_{(3, 0, 2, \cdots , 0)}(V^*)$ est de dimension $2$.

\bigskip

\bigskip

\bigskip

\bigskip

\beginsection {6. Un syst\`eme g\'en\'erateur des alg\`ebres $A(4)$ et $A(5)$.}

\bigskip

Nous pourrions, sur le mod\`ele de ce que nous avons fait concernant l'alg\`ebre $A(3)$, d\'eterminer les s\'eries de Poincar\'e de $A(4)$ et $A(5)$; mais cela ne suffirait pas pour connaitre ces alg\`ebres, aussi allons nous nous servir de certains r\'esultats de  H.W.Turnbull et J.A.Todd, d\'ej\`a \'evoqu\'es au paragraphe 2. 

\bigskip
\bigskip

\noindent {\it 6.1 Le cas de $A(4)$.}

Ici nous savons que l'alg\`ebre $A(4)$ est engendr\'ee par ses \'el\'ements de degr\'e $3$ et $6$. Par LIE, nous pouvons obtenir la structure du $\GL (4, {\Bbb C})$ module $A_6(4)$:
$$A_6(4) = \psi_{(6,0,  \cdots , 0)}(4) \oplus \psi_{(2,2,0, \cdots , 0)}(4)\oplus \psi_{(0,0,2, \cdots , 0)}(4)$$
(par exemple, nous savons que l'espace vectoriel $\psi_{(6,0,  \cdots , 0)}(4) \oplus \psi_{(2,2,0, \cdots , 0)}(4)\oplus \psi_{(0,0,2, \cdots , 0)}(4)$ est contenu dans $A_6(4)$ et ils sont tous deux de dimension $220$). Au vu de la structure de $A_6(3)$ obtenue au paragraphe 5, nous avons la:

\proclaim Proposition 6.1. L'alg\`ebre $A(4)$ est engendr\'ee par les polaris\'es de $\det x$ et de $f(\n (x), \n (y),\n (z))$.

\noindent {\bf Remarque 6.2.} A nouveau, on peut trouver par LIE, la structure du $\GL (4, {\Bbb C})$ module $A_9(4)$:

$$A_9(4) = \psi_{(9,0,  \cdots , 0)}(4) \oplus \psi_{(5,2,0, \cdots , 0)}(4)\oplus \psi_{(3,3,0, \cdots , 0)}(4)\oplus \psi_{(0,3,1, \cdots , 0)}(4)\oplus 2\psi_{(3,0,2, \cdots , 0)}(4)$$
$$\oplus \psi_{(1,0,0,2, \cdots , 0)}(4)\oplus \psi_{(2,0,1, 1,\cdots , 0)}(4)$$
$$= {\hbox{Sym}}^3({\hbox{Sym}}^3({\Bbb C}^4)) \oplus \psi_{(3,0,2, \cdots , 0)}(4)$$

\bigskip

\noindent {\it 6.2 Le cas de $A(5)$.}

L'alg\`ebre $A(5)$ est engendr\'ee par ses \'el\'ements de degr\'e $3$ et $6$ et $9$. A nouveau LIE nous donne le $\GL (5, {\Bbb C})$ module $A_6(5)$:
$$A_6(5) = \psi_{(6,0,  \cdots , 0)}(5) \oplus \psi_{(2,2,0, \cdots , 0)}(5)\oplus \psi_{(0,0,2, \cdots , 0)}(5)$$
il n'y a donc pas de nouveaux invariants de degr\'e $6$. Par contre, l'espace vectoriel $\psi_{(9,0,  \cdots , 0)}(5) \oplus \psi_{(5,2,0, \cdots , 0)}(5)\oplus \psi_{(3,3,0, \cdots , 0)}(5)\oplus \psi_{(0,3,1, \cdots , 0)}(5)\oplus 2\psi_{(3,0,2, \cdots , 0)}(5)$ est de dimension $9520$, alors que
$$\dim A_9(5) = 9570= 9520 + 50.$$
On cherche alors par LIE un candidat $\psi_{(a_1,a_2,a_3,a_4,a_5, 0, \cdots)}(5)$ pouvant intervenir dans la d\'ecomposition du $\GL (5, {\Bbb C})$ module $A_9(5)$, on trouve par \'elimination le module 
$\psi_{(0,2,0,0,1, 0, \cdots)}(5)$, c'est-\`a-dire que nous avons:
$$A_9(5) = 
\psi_{(9,0,  \cdots , 0)}(5) \oplus \psi_{(5,2,0, \cdots , 0)}(5)\oplus \psi_{(3,3,0, \cdots , 0)}(5)\oplus \psi_{(0,3,1, \cdots , 0)}(5)\oplus$$
$$ 2\psi_{(3,0,2, \cdots , 0)}(5)
\oplus  \psi_{(0,2,0,0,1, 0, \cdots)}(5)$$
$$A_9(5) =  A_3(5)A_6(5)  \oplus  \psi_{(0,2,0,0,1, 0, \cdots)}(5).$$
Todd nous donne la notation symbolique d'un invariant qui n'appartient pas \`a $A_3(5)A_6(5)$, il s'agit de:
$$[a, a', a''] [a, a', a^{(3)}] a''_{\alpha} a^{(4)}_{\alpha} a^{(4)}_{\alpha'} a^{(3)}_{\alpha'}$$
or il r\'esulte des calculs donn\'es en appendice que l'invariant 
$${\cal U}([a, a', a''] [a, a', a^{(3)}] a''_{\alpha} a^{(4)}_{\alpha} a^{(4)}_{\alpha'} a^{(3)}_{\alpha'})$$ 
est, modulo des termes appartenant \`a $A_3(5)A_6(5)$, proportionnel \`a:
$$f(\n(x) \times \n(y), (x \times z)\times (y \times t), u)$$
par cons\'equent ce dernier invariant poss\`ede  une projection non nulle sur  $\psi_{(0,2,0,0,1, 0, \cdots)}(5)$. L'alg\`ebre $A(5)$ est engendr\'ee par les polaris\'es de $\det x$, de $f(\n (x), \n (y),\n (z))$ et de $f(\n(x) \times \n(y), (x \times z)\times (y \times t), u)$. Nous pouvons r\'esumer tout ceci dans la

\proclaim Proposition 6.3. L'alg\`ebre $A(5)$ est engendr\'ee par les polyn\^omes 
$$f(x_i, x_j, x_k), \  \  f(x_i \times x_j, x_k\times x_l, x_m \times x_n)$$
et
$$ f((x_i \times x_j)\times (x_k\times x_l), (x_m \times x_n)\times (x_o\times x_p), x_q).$$

\bigskip

\bigskip

\bigskip

\bigskip

\beginsection {7. Appendices.}

\bigskip

\noindent {\it 7.1. Obtention de la relation cubique 5.2.}

\bigskip

Dans cet appendice, nous indiquons un proc\'ed\'e pour \'ecrire la relation reliant les onze invariants $f_1, \cdots f_{11}$ et cubique en $f_{11}$ de l'alg\`ebre $A(3)$ \'evoqu\'ee au paragraphe 5.

Soit $z  = (z_{ij})$ un \'el\'ement de $V$ et notons $Z_i$ les coefficients diagonaux de $\n(z)$:
$$Z_1 = z_2z_3 -z_{23}^2, \  \ Z_2 = z_1z_3 -z_{13}^2, \  \ Z_3 = z_1z_2 -z_{12}^2.$$
Nous avons $f_3 = \det z  = z_1z_2z_3 -z_1z_{23}^2 -z_2z_{13}^2 -z_3z_{12}^2 +2z_{12}z_{23}z_{13} = 
 -2z_1z_2z_3 +\sum_{i=1}^3z_iZ_i+2z_{12}z_{23}z_{13}$, relation que nous \'ecrivons sous la forme:
 $$2z_{12}z_{23}z_{13} = f_3 + 2z_1z_2z_3 -\sum_{i=1}^3z_iZ_i
= f_3 + 2\pi -V \leqno (*)$$
o\`u nous avons pos\'e $\pi = z_1z_2z_3$ et $V = \sum_{i=1}^3z_iZ_i$. Des relations donnant les coefficients $Z_i$, et en posant $\Pi = Z_1Z_2Z_3$, il vient : 
$$z_{12}^2z_{23}^2z_{13}^2 =  -\Pi + \pi^2 -\pi V + \sum_{i<j}z_iz_jZ_iZ_j.$$
La formule (*) donne alors l'\'egalit\'e: 
$$4\Pi = -V^2 + 4\sum_{i<j}z_iz_jZ_iZ_j + 2V f_3 - 4\pi \  f_3 - f^2_3.  \leqno (**)$$
que nous voyons comme une relation liant $\Pi$, $z_i$, $Z_i$ et $f_3$.

On fait alors l'observation suivante:  il y a dans la $G$-orbite d'un \'el\'ement g\'en\'erique de $3V$ un triplet de la forme $(\lambda e, \sum y_i e_i, z)$ avec $\lambda$ non nul et $y_1, y_2, y_3$ distincts. Par abus de notation, nous d\'esignerons  par $f_1, \cdots f_{11}$, les valeurs des onze invariants en ce triplet (ainsi  par exemple: $\lambda ^3 = f_1$). Il se trouve que les triplets $(z_1,z_2,z_3)$ et $(Z_1,Z_2,Z_3)$ sont enti\`erement d\'etermin\'es par ces $f_1, \cdots f_{11}$. Ainsi, on constate que $(z_1,z_2,z_3)$ est solution du syst\`eme:

$$\left\{\matrix{z_1&+&z_2&+&z_3&=&{1\over 2} \lambda^{-2} f_5\cr
y_2y_3z_1&+&y_1y_3z_2&+&y_1y_2z_3&=&{1\over 2} f_7\cr
(y_2+y_3)z_1&+&(y_1+y_3)z_2&+&(y_1+y_2)z_3&=& \lambda^{-1}f_{10}\cr
}\right.$$
qui est de Cramer, de d\'eterminant:
$\Delta = (y_1-y_3)(y_2-y_3)(y_2-y_1)$, on a d'une mani\`ere plus pr\'ecise:
$$z_1 = {1\over \Delta} (y_3-y_2) P(y_1), \ \ z_2 = {1\over \Delta} (y_1-y_3) P(y_2), \ \ z_3 = {1\over \Delta} (y_2-y_1) P(y_3)$$
avec $P(y_i) = {1\over 2} \lambda^{-2} f_5 y_i^2 - \lambda^{-1}f_{10}y_i + {1\over 2} f_7$. Cela donne pour $\pi $:
$$\pi =  -{1\over \Delta^2}P(y_1)P(y_2)P(y_3)$$
Le produit $P(y_1)P(y_2)P(y_3)$ est sym\'etrique en $y_1,y_2,y_3$ et s'exprime donc comme polyn\^ome en 
$$\sigma_1 = \sum y_i = {1\over 2} \lambda^{-2} f_4, \ \sigma_2 =  \sum_{i<j} y_iy_j = {1\over 2} \lambda^{-1}f_6, \ \sigma_3 =  y_1y_2y_3 = f_2$$
il en est de m\^eme pour $ \Delta^2$. On v\'erifie ensuite  que  $f_1^2  \Delta^2$ est un polyn\^ome en $f_1, f_2, f_4, f_6$ et que $f_1^2  P(y_1)P(y_2)P(y_3)$ est un polyn\^ome en $f_1,f_2, f_4, f_6, f_5, f_7, f_{10} $ (les v\'erifications sont un peu longues mais ne pr\'esentent pas de difficult\'es).

 \noindent De m\^eme, le triplet $(Z_1,Z_2,Z_3)$ est solution du syst\`eme:

$$\left\{\matrix{Z_1&+&Z_2&+&Z_3&=&{1\over 2} \lambda^{-1} f_8\cr
y_1Z_1&+&y_2Z_2&+&y_3Z_3&=&{1\over 2} f_9\cr
(y_2+y_3)y_1Z_1&+&(y_1+y_3)y_2Z_2&+&(y_1+y_2)y_3Z_3&=&  \lambda^{-2} f_{11}\cr
}\right.$$
de m\^eme d\'eterminant  $\Delta$; on a:$$Z_1 = {1\over \Delta} (y_3-y_2) Q(y_1, y_2y_3), \ \ Z_2 = {1\over \Delta} (y_1-y_3) Q(y_2, y_1y_3), \ \ Z_3 = {1\over \Delta} (y_2-y_1) Q(y_3, y_1y_2)$$
avec
$Q(X,Y) = {1\over 2} \lambda^{-1}f_8Y +{1\over 2}f_9 X -  \lambda^{-2}f_{11}$. Ici on trouve que: 

\noindent $\Pi = -{1\over \Delta^2}Q(y_1, y_2y_3)Q(y_2, y_1y_3)Q(y_3, y_1y_2)$, et \`a nouveau on peut constater que $f_1^2Q(y_1, y_2y_3)Q(y_2, y_1y_3)Q(y_3, y_1y_2)$ est un polyn\^ome en $f_1,f_2, f_4, f_6, f_8, f_9, f_{11}$ (de degr\'e $3$ en $f_{11}$).

\noindent Enfin, on peut faire les m\^emes manipulations sur $V$ et constater que $f_1^2  \Delta^2 V$ est un polyn\^ome en $f_1,f_2, f_4, f_6, f_5, f_7, f_{10},f_8, f_9, f_{11} $ (de  degr\'e $1$ en $f_{11}$).

Tout ceci nous conduit \`a r\'e\'ecrire la relation (**) en la multipliant par le facteur $f_1^2  \Delta^2$:
$$4 f_1^2  \Delta^2 \Pi = f_1^2  \Delta^2 (-V^2 + 4\sum_{i<j}z_iz_jZ_iZ_j )+ 2 (f_1^2  \Delta^2 V) f_3 - 4 (f_1^2  \Delta^2\pi )\  f_3 - (f_1^2  \Delta^2) f_3^2.$$
Le membre de gauche de cette \'egalit\'e est un polyn\^ome en $f_1,f_2, f_4, f_6, f_8, f_9, f_{11}$ de degr\'e $3$ en $f_{11}$, et dans le membre de droite, l'expression $2 (f_1^2  \Delta^2 V) f_3 - 4 (f_1^2  \Delta^2\pi )\  f_3 - (f_1^2  \Delta^2) f_3^2$ est un polyn\^ome en les $f_i$, de degr\'e $1$ en $f_{11}$

Il nous reste \`a \'etudier l'expression $-V^2 + 4\sum_{i<j}z_iz_jZ_iZ_j $:

$$-V^2 + 4\sum_{i<j}z_iz_jZ_iZ_j  =- \sum z_i^2 Z_i^2 + 2 \sum_{i<j}z_iz_jZ_iZ_j =$$
$$ -{1\over \Delta^4}[(y_3-y_2)^4P^2(y_1) Q^2(y_1, y_2y_3) +(y_1-y_3)^4P^2(y_2) Q^2(y_2, y_1y_3)+$$
$$(y_1-y_2)^4P^2(y_3) Q^2(y_3, y_1y_2)] 
+{2\over \Delta^4}[(y_3-y_2)^2(y_1-y_3)^2P(y_1) P(y_2)Q(y_1, y_2y_3) Q(y_2, y_1y_3)$$
$$+(y_3-y_2)^2(y_1-y_2)^2P(y_1) P(y_3)Q(y_1, y_2y_3) Q(y_3, y_1y_2) $$
$$+(y_3-y_2)^2(y_1-y_2)^2P(y_2) P(y_3)Q(y_2, y_1y_3) Q(y_3, y_1y_2)].$$
Consid\'erons la quantit\'e
$$-(y_3-y_2)^4P^2(y_1) Q^2(y_1, y_2y_3) -(y_1-y_3)^4P^2(y_2) Q^2(y_2, y_1y_3)+$$
$$2(y_3-y_2)^2(y_1-y_3)^2P(y_1) P(y_2)Q(y_1, y_2y_3) Q(y_2, y_1y_3) = $$
$$-((y_3-y_2)^2P(y_1) Q(y_1, y_2y_3) -(y_1-y_3)^2P(y_2) Q(y_2, y_1y_3))^2$$
 l'expression $(y_3-y_2)^2P(y_1) Q(y_1, y_2y_3) -(y_1-y_3)^2P(y_2) Q(y_2, y_1y_3)$ s'annulle pour $y_1 = y_2$, et par cons\'equent  le num\'erateur de $-V^2 + 4\sum_{i<j}z_iz_jZ_iZ_j $ est divisible par $(y_1-y_2)^2$, \'etant un polyn\^ome  sym\'etrique en $y_1,y_2,y_3$ il sera divisible par $\Delta^2$; en d'autres termes $ \Delta^2 (-V^2 + 4\sum_{i<j}z_iz_jZ_iZ_j )$ est polyn\^ome  sym\'etrique en $y_1,y_2,y_3$; on v\'erifie ensuite que  $f_1^2  \Delta^2 (-V^2 + 4\sum_{i<j}z_iz_jZ_iZ_j )$ est  un polyn\^ome en $f_1,f_2, f_4, f_6, f_5, f_7, f_{10},f_8, f_9, f_{11} $ de degr\'e $2$ en $f_{11}$. Tous ces calculs, quoique fastidieux ne sont pas difficiles et sont faisables par Maple.

\bigskip

\bigskip

\noindent {\it 7.2. Calcul symbolique.}

\bigskip

Ici nous indiquons bri\`evement comment montrer que l'invariant  
$${\cal U}([a, a', a''] [a, a', a^{(3)}] a''_{\alpha} a^{(4)}_{\alpha} a^{(4)}_{\alpha'} a^{(3)}_{\alpha'})$$ 
du paragraphe 6 est reli\'e \`a $ f(\n(x) \times \n(y), (x \times z)\times (y \times t), u)$. Rappelons que si  le vecteur symbolique $a = (a_1, a_2, a_3)$ repr\'esente l'\'el\'ement $x = (x_{ij})$ de $V$, l'on a:
$$x_{11} = {\cal U}(a_1^2 ),  x_{22} =  {\cal U}( a_2^2), x_{33} =  {\cal U}(a_3^2), x_{12} =  {\cal U}(a_1a_2), x_{13} =  {\cal U}(a_1a_3), x_{23} =  {\cal U}(a_2a_3).$$
$ {\cal U}$ d\'esignant l'op\'erateur ombral. 

\noindent Nous aurons besoin d'un certain nombre de lemmes; on a tout d'abord facilement le:

\proclaim Lemme 7.1.  Si $u = (u_1, u_2, u_3)$ et $v = (v_1, v_2, v_3)$ repr\'esentent  deux \'el\'ements $x$ et $y$ de $V$, alors le triplet symbolique  produit vectoriel de $u$ avec $v$ de composantes:
$$ {u_2}{v_3} - {u_3}{v_2}, \   \   \  
{u_3}{v_1} - {u_1}{v_3}, \  \  \
 {u_1}{v_2} -
{u_2}{v_1}$$
correspond \`a l'\'el\'ement $x\times y$ de l'alg\`ebre de Jordan $V$.

\noindent  En particulier, $x\times x $ et $y\times y $ sont repr\'esent\'es par:
$$\alpha = (a_2b_3 - a_3b_2,  \  a_3b_1 - a_1b_3,  \  a_1b_2 - a_2b_1)$$
et
$$\alpha' = (a'_2b'_3 - a'_3b'_2,  \  a'_3b'_1 - a'_1b'_3,  \  a'_1b'_2 - a'_2b'_1)$$
on notera  par  $(\alpha, \alpha')$ le produit vectoriel de $\alpha$ avec $\alpha'$, qui repr\'esente  donc $(x \times x)\times (y\times y)$.

\proclaim Lemme 7.2. On a la formule:
$$[({\alpha }, {\alpha '}), u, v] =  u_{\alpha } v_{\alpha '} - u_{\alpha ' } v_{\alpha }.$$

\noindent D\'eveloppons le d\'eterminant $[({\alpha }, {\alpha '}), u, v]$ suivant la derni\`ere colonne $v$:
$$[({\alpha }, {\alpha '}), u, v] = v_1(({\alpha }, {\alpha '})_2 u_3 - ({\alpha }, {\alpha '})_3u_2)$$
$$ - v_2 (({\alpha }, {\alpha '})_1 u_3 - ({\alpha }, {\alpha '})_3u_1) + v_3 (({\alpha }, {\alpha '})_1 u_2 - ({\alpha }, {\alpha '})_2u_1)$$
on ins\`ere ensuite les formules:
$$({\alpha }, {\alpha '})_1 = {\alpha }_2{\alpha '}_3 - {\alpha }_3{\alpha '}_2, ({\alpha }, {\alpha '})_2 = {\alpha }_3{\alpha '}_1 - {\alpha }_1{\alpha '}_3,
({\alpha }, {\alpha '})_3 = {\alpha }_1{\alpha '}_2 - {\alpha }_2{\alpha '}_1.$$
En mettant  $v_i{\alpha '}_i$ (resp. $v_i{\alpha }_i$) en facteur dans les termes intervenant avec le signe plus (resp. avec le signe moins), il vient:
$$[({\alpha }, {\alpha '}), u, v]   = v_1{\alpha '}_1 [u_2{\alpha }_2 + u_3{\alpha }_3] - v_1{\alpha }_1 [u_2{\alpha '}_2 + u_3{\alpha '}_3] $$
$$ + v_2{\alpha '}_2 [u_3{\alpha }_3 + u_1{\alpha }_1] - v_2{\alpha }_2 [u_3{\alpha '}_3 + u_1{\alpha '}_1] $$
$$ + v_3{\alpha '}_3 [u_2{\alpha }_2 + u_1{\alpha }_1] - v_3{\alpha }_3 [u_2{\alpha '}_2 + u_1{\alpha '}_1] .$$
Sous cette forme, on voit bien ce qu'il faut ajouter et retrancher pour avoir:
$$[({\alpha }, {\alpha '}),u, v]  = v_1{\alpha '}_1 u_{\alpha } -  v_1{\alpha }_1 u_{\alpha '}+ v_2{\alpha '}_2 u_{\alpha } - v_2{\alpha }_2 u_{\alpha '} + 
v_3{\alpha '}_2 u_{\alpha } - v_3{\alpha }_3 u_{\alpha '}$$
ce qui montre le lemme.

\noindent Si  $x$, $y$, $z$ et $t$ sont repr\'esent\'es par les vecteurs $a$(ou $b$), $a'$(ou $b'$), $a''$(ou $b''$) et $a^{(3)}$(ou $b^{(3)}$), on montre de m\^eme le: 

\proclaim Lemme 7.3. L'expression:
 $$a^{(3)} [a',a,a''] - a' [a^{(3)},a,a'']$$
repr\'esente la matrice $(x \times z)\times (y\times t)$.

\noindent Enfin, nous aurons besoin d'un dernier r\'esultat:

\proclaim Lemme 7.4. Si $A = (A_1, A_2, A_3)$, $B = (B_1, B_2, B_3)$ et $C = (C_1, C_2, C_3)$ repr\'esentent les trois \'el\'ements $X, Y, Z$ de $V$, on a:
$${\cal U}([A, B, C]^2) = f(X, Y, Z).$$

\noindent Ceci r\'esulte de la formule $f(X, Y, Z) = \tr ( (X\times Y)Z)$, du lemme 7.1, et de 
$${\cal U}(\sum a_ia'_i) = \tr(x.y)$$
o\`u $a =(a_1, a_2, a_3)$ et $a' = (a'_1, a'_2, a'_3)$ repr\'esentent $x$ et $y$.

Si maintenant $u$ est un cinqui\`eme \'el\'ement de $V$ repr\'esent\'e par $a^{(4)}$, les trois premiers lemmes nous donnent:
$$ f(\n(x) \times \n(y), (x \times z)\times (y \times t), u)= {\cal U}([(\alpha, \alpha'), a^{(3)} [a',a,a''] - a' [a^{(3)},a,a''], a^{(4)}]^2) =$$
$${\cal U}(([(\alpha, \alpha'), a^{(3)}, a^{(4)}] [a',a,a''] - [(\alpha, \alpha'), a', a^{(4)}][a^{(3)},a,a''])^2).$$
Par le lemme 7.2 nous avons:
 $$[({\alpha }, {\alpha '}), a^{(3)}, a^{(4)}] = a_{\alpha }^{(3)} a_{\alpha '}^{(4)} - a_{\alpha ' }^{(3)} a_{\alpha }^{(4)}$$
$$[({\alpha }, {\alpha '}), a', a^{(4)}] = a'_{\alpha } a_{\alpha '}^{(4)} - a'_{\alpha ' } a_{\alpha }^{(4)}.$$
On ins\`ere ces formules; en d\'eveloppant le carr\'e on constate alors facilement que l'application ombrale transforme les termes en produits d'invariants d\'ej\`a connus, sauf les deux expressions symboliques:
$$2a_{\alpha }^{(3)} a_{\alpha '}^{(4)}(a',a,a'')a'_{\alpha ' } a_{\alpha }^{(4)}(a^{(3)},a,a'')$$
et
$$2a_{\alpha ' }^{(3)} a_{\alpha }^{(4)}(a',a,a'')a'_{\alpha } a_{\alpha '}^{(4)}(a^{(3)},a,a'')$$
il est classique (voir [Gra.-You.] page 275, ou [Tod.2] page 400) que l'application ombrale transforme une expression symbolique contenant $a'_{\alpha ' } = [a', b', c']$ o\`u $a'$, $b'$ et $c'$ repr\'esentent l'\'el\'ement $y$, en un produit (\'eventuellement nul) d'invariants divisible par $\det y$. Finalement, modulo $A_3(5) A_6(5)$, nous avons le r\'esultat escompt\'e
$$f((x \times x)\times (y\times y), (x \times z)\times (y\times t),u) = 2{\cal U}(a'_{\alpha }a_{\alpha }^{(4)}a_{\alpha '}^{(4)}a_{\alpha ' }^{(3)}(a^{(3)},a,a'')(a',a,a'')).$$

\vfill

\eject

\beginsection {Bibliographie.}

\noindent [Bri.]  M. Brion, {\it Invariants et covariants des groupes r\'eductifs,} in M. Brion, G.W. Schwarz, {\it Th\'eorie des invariants et g\'eom\'etrie des vari\'et\'es quotients,} Paris, Hermann, 2000.

\noindent [Cia.]  C. Ciamberlini, {\it  Sul sistema di tre forme ternarie quadratiche,} Giornale di Battaglini 24 (1886) 141-157.

\noindent [Der.-Kem.]   H. Derksen, G. Kemper, {\it Computational Invariant Theory,} Springer-Verlag 2002.

\noindent [Far.-Kor.] J. Faraut, A. Koranyi, {\it  Analysis on Symmetric Cones,} Clarendon Press, Oxford 1994.

\noindent [Gor.1]  P. Gordan, {\it \"Uber B\"uschel von Kegelschnitten,} Math. Ann 19 (1882) 529-552.

\noindent [Gor.2]  P. Gordan, {\it Vorlesungen \"uber Invariantentheorie,} Chelsea, 1987 (reprint des originaux, Leipzig, 1885 et 1887).

\noindent [Gra.-You.]  J.H. Grace, A. Young, {\it Algebra of invariants,} Chelsea,  (reprint de l'original, Cambridge University Press, 1903).

\noindent [Gro.- Rot.- Ste.] F.D. Grosshans, G-C. Rota, J.A. Stein, {\it Invariant theory and superalgebras,} AMS Regional Confe. Series 69 (1987).

\noindent [Ilt.] A.V. Iltyakov, {\it On rational Invariants of the Group $E_6$,} Proc. Ame. Math. Soc. 124 (1996) 3637-3640.

\noindent [Kno.- Lit.] F. Knop, P. Littelmann, {\it Der Grad erzeugender Funktionen von Invariantenringen,} Mathemat. Zeit. 196 (1987) 211-229.

\noindent [Pro.] C. Procesi, {\it Lie Groups. An Approach through Invariants and Representations,} Springer-Verlag 2007.

\noindent [Rot.- Stu.]   G-C. Rota, B. Sturmfels, {\it Introduction to Invariant Theory in Superalgebras,} in D. Stanton ed., {\it Invariant Theory and Tableaux,} Springer-Verlag, IMA volumes in mathematics and its applications, vol 19, 1990.

\noindent [Sch.1]  G.W. Schwarz, {\it Representations of Simple Lie Groups with Regular Rings of Invariants,} Inventiones mathematicae 49 (1978) 167-191.

\noindent [Sch.2]  G.W. Schwarz, {\it  On Classical Invariant Theory and Binary cubics,} Ann. Inst. Fourier, 37 (1987) 191-216.

\noindent [Spr.]  T.A. Springer, {\it Jordan Algebras and Algebraic Groups,} Springer-Verlag, Ergebn. der Math. und ihrer Grenz., vol 75, 1973.

\noindent [Tod.1]    J.A. Todd, {\it The complete irreducible system of four ternary quadratic,} Proc. Lond. Math. Soc. (2) 52 (1951)   271 - 294.

\noindent [Tod.2]    J.A. Todd, {\it Ternary quadratic types,} Philos. Trans. Roy. Soc. London. Ser. A, 241 (1948)   399 - 456.

\noindent [Tur.]   H.W. Turnbull, {\it Ternary quadratic types,} Proc. London Math. Soc 2 (1910)  81-121.

\noindent [Vin.] E.B. Vinberg, {\it On Invariants of a set of matrices,}  Jour. of Lie Theory, 6, (1996) 249-269.

\noindent [Vus.] T. Vust, {\it Sur la th\'eorie des invariants des groupes classiques,} Ann. Inst. Fourier 26 (1976)  1-31.

\noindent [Wae.] B.L. Van der Waerden {\it Over het concomitantensystem van twee en drie ternaire quadratische vormen,} Amst. Ak. Versl. 32 (1923) 138-147.

\bye